
\input amstex
\documentstyle{amsppt}
\magnification = 1200
\pageheight{7.5in}
\vcorrection{-.12in}
\NoBlackBoxes
\NoRunningHeads
\nologo
\TagsOnRight
\topmatter
\title{Schwarzian Derivatives and Uniform Local Univalence}
\endtitle
\author{Martin Chuaqui, Peter Duren, and Brad Osgood}
\endauthor

\address Facultad de Matem\'aticas, P. Universidad Cat\'olica de Chile,
Casilla 306, Santiago 22, Chile
\endaddress
\email mchuaqui\@mat.puc.cl
\endemail

\address Department of Mathematics, University of Michigan, Ann Arbor,
Michigan 48109--1043
\endaddress
\email duren\@umich.edu
\endemail
 
\address Department of Electrical Engineering, Stanford University,
Stanford, California 94305
\endaddress
\email osgood\@ee.stanford.edu
\endemail

\dedicatory{Dedicated to Walter Hayman on the occasion of his 80{\it th} birthday}
\enddedicatory

\thanks The authors are supported by Fondecyt Grant \# 1030589. 
\endthanks

\abstract  Quantitative estimates are obtained for the (finite) valence of functions 
analytic in the unit disk with Schwarzian derivative that is bounded or of slow 
growth.  A harmonic mapping is shown to be uniformly locally univalent with 
respect to the hyperbolic metric if and only if it has finite Schwarzian norm, 
thus generalizing a result of B. Schwarz for analytic functions.  A numerical 
bound is obtained for the Schwarzian norms of univalent harmonic mappings. 
\endabstract

\subjclass Primary 30C99, Secondary 31A05, 30C55
\endsubjclass
\keywords Analytic function, valence, harmonic mapping, Schwarzian derivative, uniform 
local univalence, Schwarzian norm, minimal surface, harmonic lift 
\endkeywords

\endtopmatter

\document

\flushpar
{\bf \S 1.  Finite valence.}
\smallpagebreak

     Our point of departure is a classical theorem of Nehari [14] that gives a general criterion for univalence of an analytic function in terms of its Schwarzian derivative
$$
{\Cal S}f = (f''/f')' - \tfrac12(f''/f')^2\,. 
$$
A positive continuous even function $p(x)$ on the interval $(-1,1)$ is called a {\it Nehari function} if 
$(1-x^2)^2p(x)$ is nonincreasing on $[0,1)$ and no nontrivial solution $u$ of the differential 
equation $u''+pu=0$ has more than one zero in $(-1,1)$.   Nehari's theorem can be stated as 
follows.
\proclaim{Theorem A}   Let $f$ be analytic and locally univalent in the unit disk \,$\Bbb D$, and 
suppose its Schwarzian derivative satisfies 
$$
|{\Cal S}f(z)|\leq 2p(|z|)\,, \qquad z\in{\Bbb D}\,,  \tag1
$$
for some Nehari function $p(x)$.  Then $f$ is univalent in $\Bbb D$.  
\endproclaim

     As special cases the theorem includes the criteria  $|{\Cal S}f(z)|\leq2(1-|z|^2)^{-2}$ 
and $|{\Cal S}f(z)|\leq{\pi}^2/2$ obtained earlier by Nehari [13], as well as the criterion 
$|{\Cal S}f(z)|\leq4(1-|z|^2)^{-1}$ stated by Pokornyi [17].  The weaker inequality 
$$
|{\Cal S}f(z)|\leq \frac{2(1+\delta^2)}{(1-|z|^2)^2}\,, 
\qquad z \in \Bbb D\,, 
$$
does not imply univalence, but it does imply uniform local univalence in the sense that 
the hyperbolic distance $d(\alpha,\beta)\geq\pi/\delta$ for any pair of points 
$\alpha,\beta\in\Bbb D$ where $f(\alpha)=f(\beta)$.   In  a previous paper [4] we gave a 
streamlined proof of this result, which is due to B. Schwarz [19], and demonstrated the sharpness 
of the lower bound (see also Minda [12]).   Furthermore, we showed that any weaker form 
$|{\Cal S}f(z)|\leq C\,p(|z|)$ of Nehari's condition (1) still implies that $f$ has finite  
valence if $(1-x^2)^2p(x)\to0$ as $x\to1-$.  In particular, if $|{\Cal S}f(z)|\leq C$ for some  
constant $C$ and all $z\in\Bbb D$, then $f$ has finite valence in the unit disk.  

     We now derive this last result by a more elegant method, which also provides a quantitative 
bound for the valence in terms of the constant $C$.   By the {\it valence} of $f$ we mean 
$N=\sup_{w\in\Bbb C} n(f,w)$, where $n(f,w)\leq\infty$ is the number of points $z\in\Bbb D$ for 
which $f(z)=w$.   Here is our theorem.  
\proclaim{Theorem 1}  Let $f$ be analytic and locally univalent in the unit disk $\Bbb D$, and 
suppose its Schwarzian derivative satisfies 
$$
|{\Cal S}f(z)|\leq C\,, \qquad z\in{\Bbb D}\,,  
$$
for some constant $C>\pi^2/2$\,.  Then $|\alpha-\beta|\geq\sqrt{2/C}\,\pi$ for any pair of points 
$\alpha,\beta\in\Bbb D$ where $f(\alpha)=f(\beta)$.  Consequently, $f$ has finite 
valence and assumes any given value at most $\left(1+\frac{\sqrt{2C}}{\pi}\right)^2$ times. 
\endproclaim

     Before embarking on the proof, we recall some standard facts about the Schwarzian 
derivative.  It is M\"obius invariant: ${\Cal S}(T\circ f)={\Cal S}f$ for every M\"obius transformation 
$$
T(z) = \frac{az+b}{cz+d}\,,\qquad ad-bc\neq 0\,.
$$
Also, ${\Cal S}(f\circ T)=(({\Cal S}f)\circ T){T'}^2$\,.
For any analytic function $\psi$, the functions $f$ with Schwarzian ${\Cal S}f=2\psi$ are 
precisely those of the form $f=u_1/u_2$, where $u_1$ and $u_2$ are linearly independent 
solutions of the differential equation $u''+ \psi u=0$.  Thus if ${\Cal S}f=2\psi$, then 
$f(\alpha)=f(\beta)$ if and only if some solution of the differential equation $u''+ \psi u=0$ 
vanishes at $\alpha$ and $\beta$.

     We will make use of the following lemma, which is a variant of a lemma in [7].  
\proclaim{Lemma 1}  Suppose that $u=u(z)$ is a solution of the differential equation 
$u''+ \psi u=0$  for some function $\psi$ analytic in $\Bbb D$.  Let $z=z(s)\,,\ s\in (0,b)$\,, be 
an arclength parametrization of a line segment in $\Bbb D$, and suppose that 
$v(z)=|u(z(s))|>0$ for $s$ in the interval $(0,b)$.  Then 
$$
v''(s) + |\psi(z(s))|\,v(s) \geq 0\,, \qquad 0<s<b\,.
$$
\endproclaim
\demo{Proof of lemma}  Differentiation of $v^2=u\overline{u}$ gives 
$$
v(s)\,v'(s) = \text{Re}\bigl\{u'(z(s))\,z'(s)\,\overline{u(z(s))}\,\bigr\}\,, \qquad 0<s<b\,.
$$
But $v(s)>0$ and $|z'(s)|=1$, so it follows that 
$$
v(s)\,|v'(s)| \leq |u'(z(s))|\,v(s)\,, \qquad \text{or} \qquad  |v'(s)| \leq |u'(z(s))|\,.
$$
Differentiation of $vv'$ gives
$$ 
vv'' + {v'}^2 = \text{Re}\bigl\{u''{z'}^2\,\overline{u} + |u'|^2\bigr\}\,,
$$
since $|z'(s)|=1$ and $z'(s)$ is constant for the parametrization of a line segment.  
Introducing the differential equation $u''=-\psi u$, we conclude that 
$$
\align
vv'' + {v'}^2 &= |u'|^2 - \text{Re}\bigl\{ \psi\, |u|^2 {z'}^2\bigr\} \geq  |u'|^2 - |\psi| |u|^2 \\
&\geq |v'|^2 - |\psi| |u|^2 = {v'}^2 - |\psi| v^2 \,. 
\endalign
$$
Therefore, $v(v'' + |\psi| v) \geq  0$\,, and the desired result follows because $v(s)>0$ on the 
interval $(0,b)$.
\qed\enddemo
\demo{Proof of theorem}  Under the hypothesis $|\psi(z)|\leq C/2$, where ${\Cal S}f=2\psi$, suppose 
that $f(\alpha)=f(\beta)$ for some pair of distinct points $\alpha, \beta \in \Bbb D$.  Then some 
solution of the differential equation $u''+\psi u=0$ vanishes at $\alpha$ and $\beta$.  Without loss 
of generality, we may suppose that $u(z)\neq0$ on the open line segment with endpoints 
$\alpha$ and $\beta$.  Let $z=z(s)$ be the parametrization of this segment by arclength $s$, 
where $z(0)=\alpha$ and $z(b)=\beta$, so that $b=|\alpha-\beta|$.  Then by Lemma 1, the 
function $v(s)=|u(z(s))|$ has the properties $v(0)=v(b)=0$, $v(s)>0$, and 
$$
v''(s) + |\psi(z(s))|\,v(s) \geq 0\,, \qquad 0 < s < b\,.
$$
We now apply the Sturm comparison theorem (see for instance [1]).  Note that $v(s)$ is a 
real-valued function that satisfies the differential equation $v''(s)+g(s)v(s)=0$, with
$$
g(s) = - v''(s)/v(s) \leq |\psi(z(s))| \leq C/2\,.
$$
On the other hand, the solutions of the differential equation $y'' + (C/2)y=0$ are sinusoids whose 
zeros are separated by the distance $\sqrt{2/C}\,\pi$.  By the Sturm comparison theorem, 
$$
|\alpha-\beta| = b \geq \sqrt{2/C}\,\pi\,,
$$
as claimed.  Note that if $C = \pi^2/2$, then the argument shows that $|\alpha-\beta| \geq2$, 
and so we recover Nehari's theorem that $f$ is univalent in $\Bbb D$ if $|{\Cal S}f(z)|\leq\pi^2/2$\,.

     Now for the estimate of valence.  Let $w$ be an arbitrary complex number.  By what we 
have already proved, the points in $\Bbb D$ where $f(z)=w$ are the centers of disjoint 
disks of radius $\pi/\sqrt{2C}$\,.  If there are $N$ such points, a comparison of areas shows that 
$$
N\,\pi\left(\frac{\pi}{\sqrt{2C}}\right)^2 \leq \pi\left(1 + \frac{\pi}{\sqrt{2C}}\right)^2\,,
$$
which reduces to the stated inequality $N\leq \left(1+\frac{\sqrt{2C}}{\pi}\right)^2$\,.  
\qed\enddemo

     The bound on the valence is not sharp.  For instance, for $C=\pi^2/2$ it gives $n\leq4$, 
whereas Nehari's theorem shows that $n\leq1$.   Nevertheless, the question remains whether the 
bound is sharp in order of magnitude.  Theorem 1 shows that under the condition $|{\Cal S}f(z)|\leq C$ 
the sharp bound on the valence is $O(C)$ as $C\to\infty$.  On the other hand, the simple example 
$$
f(z) = \tan\left(\sqrt{C/2}\,z\right)\,, \qquad \text{for which} \quad {\Cal S}f(z) = C > \pi^2/2\,, \tag2
$$
shows that the valence may increase as fast as $\sqrt{C}$.  Indeed, $f(x)=0$ for all points 
$x=\pm k\pi\sqrt{2/C}$ where $k=0,1, 2, \dots$, and at least $\frac{\sqrt{2C}}{\pi} -1$ of 
these points lie in the unit disk.  Thus the bound on the valence cannot be improved to 
anything better than $O(\sqrt{C})$ as $C\to\infty$.

\bigpagebreak
\flushpar
{\bf \S 2.  Schwarzians of slow growth.}   
\smallpagebreak

      We showed in [4] that for each Nehari function $p(x)$ with $(1-x^2)^2p(x)\to0$ as $x\to1-$, 
any condition of the form $|{\Cal S}f(z)|\leq C\,p(|z|)$ implies that $f$ has finite valence in the disk.  
In the previous section we considered functions with $|{\Cal S}f(z)|\leq C$ and obtained  an explicit estimate for the valence in terms of $C$.   We now take $p(x)= \frac{2}{1-x^2}$, the Nehari function 
in the univalence criterion of Pokornyi [17], and derive an estimate, in terms of the constant $C$, 
for the (finite) valence of functions $f$ with $|{\Cal S}f(z)|\leq C\,p(|z|)$.  We will content ourselves 
with an asymptotic estimate as $C\to\infty$, although the proof can be adapted to yield an explicit 
bound.  

\proclaim{Theorem 2}   Let $f$ be analytic and locally univalent in $\Bbb D$, and suppose its 
Schwarzian derivative satisfies
$$
|{\Cal S}f(z)|\leq \frac{2C}{1-|z|^2} \,, \qquad z\in\Bbb D\,, \tag3
$$
for a constant $C>2$.  Then $f$ has finite valence $N = N(C) \leq A\,C \log C$, where $A$ is some 
absolute constant.     
\endproclaim

     The proof of Theorem 2 will invoke the separation result of Theorem 1.  The following 
geometric lemma will be useful.  

\proclaim{Lemma 2}  If \,$n$ points $z_1, z_2, \dots, z_n$ lie in an annulus 
$\rho\leq|z|\leq\rho+d\leq1$ and have the separation property $|z_j - z_k|\geq 2d$ for $j\neq k$, 
then $n\leq2\pi/d$.  
\endproclaim
\demo{Proof of lemma}  It will suffice to show that $|\arg\{z_j\} - \arg\{z_k\}| > d$ for $j\neq k$.  
But if $|\arg\{z_j\} - \arg\{z_k\}| \leq d$ for some $j\neq k$, then by the triangle inequality 
$$
|z_j - z_k| \leq d + \rho d < 2d\,, 
$$
which contradicts the hypothesis.
\qed\enddemo
\demo{Proof of theorem}  In terms of the Nehari function  $p(x)= 2/(1-x^2)$\,, 
the hypothesis is that $|{\Cal S}f(z)|\leq C\,p(|z|)$.  We claim that $f$ is univalent in the disk 
$|z|<r_0=\pi/\sqrt{{\pi}^2+ 4C}$.  Indeed, the function $g(z)=f(r_0z)$ has Schwarzian 
${\Cal S}g(z)=r_0^2\,{\Cal S}f(r_0z)$, and so 
$$
|{\Cal S}g(z)| \leq r_0^2\, C\, p(r_0) = \frac{{\pi}^2}{2}\,,
$$
which implies that $g$ is univalent in $\Bbb D$, by Nehari's theorem.  Thus $f$ is univalent 
in $|z|<r_0$.

     We now define the sequence $\{r_k\}$ recursively by the formula  
$$
r_k - r_{k-1} = d_k = \frac{\pi}{\sqrt{2C\,p(r_k)}}\,, \qquad k=1,2,\dots\,. \tag4
$$
If $r_k<1$, then since 
$|{\Cal S}f(z)|\leq C\,p(r_k)$ in the disk $|z|\leq r_k$, the Schwarzian of $g(z)=f(r_kz)$ satisfies 
$|{\Cal S}g(z)|\leq C\,r_k^2\,p(r_k)$ in $\Bbb D$.  Thus by Theorem 1, if $f(\alpha)=f(\beta)$ for 
two points $\alpha$ and $\beta$ in the disk $|z|<r_k$, then 
$$
|\alpha - \beta| \geq \frac{r_k\,\sqrt{2}\,\pi}{\sqrt{C\,r_k^2\,p(r_k)}} =  
\frac{\sqrt{2}\,\pi}{\sqrt{C\,p(r_k)}} = 2\,d_k\,.
$$
An appeal to Lemma 2 now shows that the valence $N_k$ of $f$ in the annulus 
$r_{k-1}\leq|z|<r_k$ satisfies 
$$
N_k \leq \frac{2\pi}{d_k} = 2 \sqrt{2C\,p(r_k)}\,. \tag5
$$

     Next we make a closer examination of the recurrence relation (4), which we rewrite as 
$$
x - a = \varepsilon  \, \sqrt{1 - x^2}\,, \qquad \text{where} \ \ a = r_{k-1}\,, \ \ x=r_k\,, \ \ \text{and} \ \ 
\varepsilon = \frac{\pi}{2\sqrt{C}}\,.
$$
Squaring and solving the quadratic equation, we find 
$$
x = \frac{1}{1+\varepsilon^2} \left(a + \varepsilon \sqrt{1 - a^2 + \varepsilon^2}\right)\,,
$$
which leads after further calculation to the formula 
$$
\frac{1}{x-a} = \frac{\sqrt{1 - a^2 + \varepsilon^2} + \varepsilon a}{\varepsilon(1 - a^2)} 
= \phi(a)\,, \quad \text{say}.
$$
It is important to observe that $\phi$ is an increasing function on the interval $0<a<1$.   
This be verified by computing its derivative:
$$
\varepsilon(1-a^2)^2\sqrt{1 - a^2 + \varepsilon^2}\,{ \phi}'(a) = a(1-a^2) + 2a\varepsilon^2  
+ \varepsilon(1+a^2)\sqrt{1 - a^2 + \varepsilon^2} > 0\,.
$$
Reverting to the original notation, we have $d_k\phi(r_{k-1})=1$.  
  
       Now let $R = \left(1 - \frac{1}{4C}\right)^{1/2}$, and observe that $R>R_0$ because 
$C>2$.  Define the index $m$ by the condition $r_m\leq R<r_{m+1}$.  By virtue of (5), the 
valence of the function $f$ in the disk $|z|\leq r_m$  is bounded by 
$$
\align
1 + \sum_{k=1}^m N_k &\leq 1 + 2\pi \sum_{k=1}^m \frac{1}{d_k} = 1 + 2\pi \sum_{k=1}^m 
\phi(r_{k-1})^2 (r_k - r_{k-1}) \\
&\leq 1 + 2\pi \int_0^R \phi(x)^2\,dx\,, 
\endalign
$$
since $\phi$ is an increasing function.  Thus we need to estimate the integral 
$$
\align
\int_0^R  \phi(x)^2\,dx  &=   \frac{1}{\varepsilon^2} \int_0^R \frac{(\sqrt{1+ \varepsilon^2 - x^2} 
+ \varepsilon x)^2}{(1-x^2)^2} \, dx \\
&= \frac{1}{\varepsilon^2} \int_0^R \frac{dx}{1 - x^2} + \int_0^R \frac{1 + x^2}{(1 - x^2)^2}\,dx +  
 \frac{2}{\varepsilon} \int_0^R \frac{x\,\sqrt{1+ \varepsilon^2 - x^2}}{(1-x^2)^2} \, dx \\  
 \vspace{1 \jot}
 &= I_1 + I_2 + I_3\,, \qquad \text{say}.
\endalign
$$
 
      Recall that $R^2 = 1 - \frac{1}{4C}$ and $\varepsilon = \frac{\pi}{2\sqrt{C}}$\,, so that 
$$
\align 
I_1 &\leq \frac{4C}{\pi^2} \int_0^R \frac{dx}{1 - x^2} = \frac{2C}{\pi^2} \log\frac{(1+R)^2}{1-R^2} \\
&\leq \frac{2 \,C}{\pi^2}\,\log (16\,C) = O(C \log C)\,.
\endalign
$$
On the other hand, 
$$
I_2 = O\left(\frac{1}{1 - R^2}\right) = O(C)\,,
$$
whereas an integration by parts gives 
$$
\align
I_3 &= \frac{1}{\varepsilon} \left\{\left[\frac{\sqrt{1+ \varepsilon^2 - x^2}}{1 - x^2}\right]_0^R 
+ \int_0^R \frac{x\,dx}{(1-x^2)\sqrt{1+ \varepsilon^2 - x^2}}\right\} \\
&\leq O(C) +  \frac{1}{\varepsilon} \int_0^R \frac{x\,dx}{(1-x^2)^{3/2}} = O(C)\,.
\endalign
$$
If $r_m < R$, the same argument that produced the estimate (5) shows that in the annulus 
$r_m \leq |z| < R$ the valence of $f$ is no greater than $2\sqrt{2C\,p(R)} = O(C)$. 

     To complete the proof, we need to estimate the valence of $f$ in the annulus $R\leq|z|<1$.  
The radius $R = \left(1 - \frac{1}{4C}\right)^{1/2}$ is chosen so that $(1-R^2)^2C\,p(R)=\frac12$.   
The radius $R_1 = \left(1 - \frac{1}{2C}\right)^{1/2}$ has the properties $0<R_1<R$ and 
$$
(1-R_1^2)^2C\,p(R_1) = 2C(1-R_1^2) = 1 \,.
$$
Thus the bound (3) on the Schwarzian derivative of $f$ implies that 
$$
|{\Cal S}f(z)| \leq \frac{1}{(1-|z|^2)^2} \,, \qquad R_1 \leq |z| < 1\,. \tag6 
$$
Suppose now that $f(\alpha)=f(\beta)$ for two points $\alpha$ and $\beta$ in the annulus 
$R\leq|z|<1$.  Then by Nehari's theorem, or rather by its proof, the hyperbolic geodesic joining $\alpha$ and $\beta$ cannot lie entirely in the annulus $R_1\leq|z|<1$.  For then the Schwarzian of $f$ would satisfy (6) along such a geodesic. By a well-known technique of Nehari [14], this would imply that a function $g=f\circ\varphi$, where $\varphi$ is a suitable conformal automorphism of the disk, satisfies 
$|{\Cal S}g(x)|\leq (1-x^2)^{-2}$ on the real interval $-1<x<1$ and has the property $g(a)=g(b)$ for a 
pair of distinct points $a$ and $b$ in that interval.  Equivalently, a solution to the associated linear 
differential equation vanishes at two points of the interval $(-1,1)$, which is not possible.  This shows 
that  $f$ is univalent in each part of the annulus $R\leq|z|<1$ which lies inside the arch of some hyperbolic geodesic entirely contained in the larger annulus $R_1\leq|z|<1$.  The conclusion 
is strongest if we take the hyperbolic geodesic to be tangent to the circle $|z|=R_1$.  
     
     The estimate of valence in the annulus $R\leq|z|<1$ now reduces to a covering problem, namely 
to estimate the number of curvilinear rectangles required to cover the annulus.  Here a   
{\it curvilinear rectangle} is understood to mean the intersection of the given annulus with the 
region inside a hyperbolic geodesic that is tangent to the circle $|z|=R_1$.  Observe that the geodesic 
that is tangent to this circle at the point $z=R_1$ is the image of the imaginary axis under the M\"obius 
automorphism 
$$
 T(z) = \frac{z+R_1}{1+R_1z}\,, \qquad z\in\Bbb D\,.
$$
In order to locate the two points where this geodesic meets the circle $|z|=R$\,, we calculate that 
$|T(iy)|=R$ implies 
$$
y^2 = \frac{R^2 - R_1^2}{1-R^2R_1^2} = \frac{2C}{6C-1}\,.
$$
Choosing $y>0$, we find by further calculation that 
$$
\align
\arg\{T(iy)\} &= \tan^{-1}\left(\frac{y}{1+y^2} \, \frac{1-R_1^2}{R_1}\right) \\
&= \tan^{-1}\left(\frac{\sqrt{6C-1}}{8C-1} \, \frac{1}{\sqrt{2C-1}}\right) \geq 
\tan^{-1}\left(\frac{\sqrt{3}}{8C}\right) \geq \frac{1}{5C}
\endalign
$$
for all constants $C$ sufficiently large.  Therefore, the annulus $R\leq|z|<1$ is contained in 
the union of at most $[5\pi C] +1$ curvilinear rectangles of the type described, where $[x]$ 
denotes the integer part of $x$.  Consequently, the valence of $f$ in this annulus is $O(C)$ 
as $C\to\infty$.   This concludes the proof of Theorem 2.
\qed\enddemo

     The example (2) again shows that the estimate of valence in Theorem 2 cannot be improved   
to $o(\sqrt{C})$.  In search of a better lower bound, it is natural to investigate  the zeros of solutions 
of the differential equation 
$$
y'' + \frac{C}{1-x^2} \,y = 0 \tag7
$$
in the interval $(-1,1)$.  The solutions of (7) are easily seen to have the form $y=(1-x^2)u'$, 
where $u$ is a solution of the Legendre equation 
$$
(x^2 - 1) u''(x) + 2x\,u'(x) - C\,u(x) = 0\,. \tag8
$$
(Compare Kamke [10], 2.240, eq\. 14, p\. 460.)  If $C=n(n+1)$ for $n=1,2,\dots$, one solution 
of (8) is the Legendre polynomial $u=P_n(x)$, which is known to have exactly $n$ simple 
zeros in the interval $(-1,1)$.  Thus by Rolle's theorem, the derivative $P_n'(x)$, a polynomial of 
degree $n-1$, has exactly $n-1$ zeros in $(-1,1)$.  In other words, if $C=n(n+1)$, then some 
solution of (7) has at least $n-1$ zeros in the unit disk.  This remains true, by the Sturm comparison 
theorem, if $n(n+1)<C<(n+1)(n+2)$.  The conclusion is that some analytic function whose 
Schwarzian satisfies (3) has valence (loosely speaking) at least $\sqrt{C}$, which shows again 
that the asymptotic estimate of Theorem 2 cannot be improved beyond $N=O(\sqrt{C})$ as 
$C\to\infty$.  

\bigpagebreak
\flushpar
{\bf \S 3.  Uniform local univalence and harmonic mappings.}   
\smallpagebreak

     The {\it pseudohyperbolic metric} is defined by  
$$
\rho(\alpha,\beta)=\left|\frac{\alpha-\beta}{1-\overline{\alpha}\beta}
\right|\,, \qquad \alpha, \beta\in{\Bbb D}\,, 
$$     
and is M\"obius invariant.  More precisely, $\rho(\varphi(\alpha),\varphi(\beta)) = \rho(\alpha,\beta)$ 
if $\varphi$ is any M\"obius self-mapping of $\Bbb D$.  The {\it pseudohyperbolic disk} with 
center $\alpha$ and radius $r$ is defined by 
$$
\Delta(\alpha,r) = \bigl\{ z\in\Bbb D : \rho(z,\alpha) < r \bigr\}\,.
$$
It is a true Euclidean disk, but $\alpha$ and $r$ are not the Euclidean center and radius unless 
$\alpha=0$.  The {\it hyperbolic metric} is 
$$
d(\alpha, \beta) = \frac12 \, \log\frac{1 + \rho(\alpha, \beta)}
{1 - \rho(\alpha, \beta)}\,.
$$

     The {\it Schwarzian norm} of a function $f$ analytic and locally univalent in the unit disk  
is defined by
$$
\|{\Cal S}f\| = \sup_{z\in\Bbb D} \,(1-|z|^2)^2 |{\Cal S}f(z)|\,.
$$
It is M\"obius invariant in the sense that $\|{\Cal S}(f\circ\varphi)\| = \|{\Cal S}f\|$ for any 
M\"obius self-mapping $\varphi$ of the unit disk.  The previously mentioned result of 
Nehari [13], a special case of Theorem A, can be rephrased to say that $f$ is univalent 
in $\Bbb D$ if $\|{\Cal S}f\| \leq2$.  In the converse direction, Kraus [11] 
showed that $\|{\Cal S}f\|\leq6$ whenever $f$ is analytic and univalent in $\Bbb D$. 
The bound is sharp, since the Koebe function $k(z)=z/(1-z)^2$ has Schwarzian 
$$
{\Cal S}k(z) = - \,\frac{6}{(1-z^2)^2}\,.
$$

     According to the theorem of B. Schwarz [19], the condition $\|{\Cal S}f\|<\infty$ implies that 
$f$ is {\it uniformly locally univalent} in the hyperbolic metric, or equivalently in the pseudohyperbolic 
metric.  Specifically, this means that for some radius $r>0$, the function $f$ is univalent in every 
pseudohyperbolic disk $\Delta(\alpha,r)$.  Conversely, if $f$ is uniformly locally univalent, 
then $\|{\Cal S}f\|<\infty$.  In fact, it is known that $\|{\Cal S}f\|\leq6/r^2$.  To see this, 
suppose that $f$ is univalent in every pseudohyperbolic disk $\Delta(\alpha,r)$.   For any fixed 
$\alpha\in\Bbb D$, the M\"obius transformation 
$$
\varphi(z) = \frac{rz + \alpha}{1 + \overline{\alpha}rz}\,, \qquad z\in\Bbb D\,,
$$
maps $\Bbb D$ onto $\Delta(\alpha,r)$.  Thus $g=f\circ\varphi$ is univalent in $\Bbb D$, and so 
$\|{\Cal S}g\|\leq6$ by Kraus' theorem.  In particular, $|{\Cal S}g(0)|\leq6$.  But 
$$
{\Cal S}g(0) = {\Cal S}(f\circ\varphi)(0) = (({\Cal S}f)(\varphi(0)))\varphi'(0)^2 
= r^2(1-|\alpha|^2)^2 {\Cal S}f(\alpha)\,,
$$
and so $r^2(1-|\alpha|^2)^2 |{\Cal S}f(\alpha)|\leq6$.  Taking the supremum over all 
$\alpha\in{\Bbb D}$, we conclude that $\|{\Cal S}f\|\leq6/r^2$.

     To what extent do these relations generalize to harmonic mappings?  A complex-valued 
harmonic function in a simply connected domain has the canonical representation 
$f=h+\overline{g}$, unique up to an additive constant, where $h$ and $g$ are analytic 
functions.  By a theorem of H. Lewy (see [9]), the Jacobian $|h'|^2-|g'|^2$ of a locally 
univalent harmonic mapping never vanishes.  The harmonic mappings with positive 
Jacobian are said to be {\it orientation-preserving}.  These are harmonic mappings whose 
{\it dilatation} $\omega =g'/h'$ is an analytic function with $|\omega(z)|<1$.  
An orientation-preserving harmonic mapping lifts to a mapping $\widetilde{f}$ onto a minimal 
surface described by conformal parameters, if and only if $\omega=q^2$, the square of some 
analytic function $q$.  For such mappings $f$ we have defined [2]  the {\it Schwarzian 
derivative} by the formula 
$$
{\Cal S}f = 2\bigl(\sigma_{zz} - \sigma_z^2\bigr)\,,
$$
where $\sigma = \log\bigl(|h'|+|g'|\bigr)$ and  
$$
\sigma_z = \frac{\partial\sigma}{\partial z}  
= \frac12 \left(\frac{\partial\sigma}{\partial x} 
- i \frac{\partial\sigma}{\partial y}\right)\,, \qquad z = x+iy\,.
$$
If $f$ is analytic, ${\Cal S}f$ is the classical Schwarzian.  If $f$ 
is harmonic and $\varphi$ is analytic, then $f\circ\varphi$ is harmonic and 
$$
{\Cal S}(f\circ\varphi) = (({\Cal S}f)\circ\varphi){\varphi'}^2 + {\Cal S}\varphi\,,
$$
a generalization of the classical formula for analytic functions $f$.  In particular, 
$$
{\Cal S}(f\circ\varphi) = (({\Cal S}f)\circ\varphi){\varphi'}^2 
$$
if $\varphi$ is a M\"obius self-mapping of the disk.  From this it follows that the 
Schwarzian norm 
$$
\|{\Cal S}f\| = \sup_{z\in\Bbb D} \,(1-|z|^2)^2 |{\Cal S}f(z)|\,.
$$
of a harmonic mapping retains the M\"obius invariance property 
$\|{\Cal S}(f\circ\varphi)\| = \|{\Cal S}f\|$.
\proclaim{Theorem 3}  Let $f=h+\overline{g}$ be an  orientation-preserving harmonic 
mapping whose dilatation $\omega =g'/h'$ is the square of an analytic function in the unit 
disk.  Then $\|{\Cal S}f\|<\infty$ if and only if $f$ is uniformly locally univalent.  
\endproclaim

     The proof will invoke a recent result of Chuaqui and Hern\'andez [6], which we state here 
for reference.  
\proclaim{Theorem B}   Let $f=h+\overline{g}$ be an  orientation-preserving harmonic 
mapping in the unit disk, and suppose that $h$ is univalent and $h(\Bbb D)$ is convex.  
Then $f$ is univalent in $\Bbb D$.  
\endproclaim
\demo{Proof of Theorem B}  In the paper [6] this result comes out of a more general argument, 
but the proof for this special case is so short that we include it here for completeness.  
If $f(z_1)=f(z_2)$, then $h(z_1)-h(z_2)=\overline{g(z_2)}-\overline{g(z_1)}$.  With the notation 
$w_1=h(z_1)$ and $w_2=h(z_2)$, this can be written as 
$$
\overline{w_1} -  \overline{w_2} = \varphi(w_2) - \varphi(w_1)\,, \qquad \text{where} \ \ 
\varphi=g\circ h^{-1}\,.
$$
But $\varphi$ is analytic on the convex domain $h(\Bbb D)$, so this says that 
$$
\overline{w_1} -  \overline{w_2} = \int_{w_1}^{w_2} \varphi'(w)\,dw\,,
$$
where the integral is taken over a straight-line segment.  However, this is not possible, because 
$|\varphi'(w)|=|g'(z)/h'(z)|<1$ by the hypothesis that $f$ is orientation-preserving.  
\qed\enddemo

     We will also need a result that is implicit in work of Sheil-Small [20].  An analytic or harmonic 
function $f$ is said to be {\it uniformly locally convex} if there exists a radius $r>0$ such that 
$f$ maps every pseudohyperbolic disk $\Delta(\alpha,r)$ univalently onto a convex region.  
\proclaim{Theorem C}  Let $f=h+\overline{g}$ be an  orientation-preserving harmonic 
mapping that is uniformly locally univalent in the unit disk.  Then its analytic part $h$ is 
uniformly locally convex.  
\endproclaim
\demo{Proof of Theorem C}  Suppose first that $f$ is univalent in the entire disk $\Bbb D$.  Then 
we may assume without loss of generality that $f\in S_H$, the class of orientation-preserving 
univalent harmonic mappings of $\Bbb D$ for which $h(0)=g(0)=0$ and $h'(0)=1$.  The analytic 
part of such a mapping has the power series expansion $h(z)=z+a_2z^2+\dots$, and it is a 
result of Clunie and Sheil-Small that the coefficients $a_2$ have an absolute bound; in other 
words, $\lambda=\sup_{f\in S_H}  |a_2|$ is finite.  It is conjectured that $\lambda=3$, but the 
best bound currently known (see [9], p. 97) is approximately $\lambda<49$.  Now if $f\in S_H$, 
then for each fixed $\zeta\in\Bbb D$ the function 
$$
F(z) = \frac{f\left(\frac{z+\zeta}{1+\overline{\zeta}z}\right) - f(\zeta)}{(1 - |\zeta|^2)h'(\zeta)} 
= H(z)+\overline{G(z)}
$$
also belongs to the class $S_H$\,, so that $|\frac12 H''(0)|\leq\lambda$.  But a calculation gives 
$$
H''(0) = (1 - |\zeta|^2)\frac{h''(\zeta)}{h'(\zeta)} - 2\overline{\zeta}\,,
$$
so we have 
$$
\left|\frac{\zeta h''(\zeta)}{h'(\zeta)}  - \frac{2\rho^2}{1-\rho^2}\right| \leq \frac{2\lambda \rho}{1-\rho^2}\,,
$$
which implies that 
$$
\text{Re}\left\{1 + \frac{\zeta h''(\zeta)}{h'(\zeta)}\right\} \geq  \frac{1 - 2\lambda \rho + \rho^2}
{1 - \rho^2} > 0 
$$
for $|\zeta| = \rho <\mu=\lambda - \sqrt{\lambda^2 - 1}$.   By the familiar analytic criterion for 
convexity (see [8], p. 42), this shows that $h(z)$ is convex in the disk $|z|<\mu$.   If $f$ is 
assumed to be univalent only in the subdisk $|z|<r$, the preceding result can be adapted to 
show that $h$ is univalent in the disk $|z|<\mu r$.  If $f$ is univalent in the pseudohyperbolic 
disk $\Delta(\alpha,r)$, then for a suitable M\"obius self-mapping $\varphi$ of $\Bbb D$ the 
composite function $\Phi=f\circ\varphi$ is univalent in $\Delta(0,r)$, and so its analytic part is 
convex in $\Delta(0,\mu r)$, which implies that the analytic part $h$ of $f=\Phi\circ\varphi^{-1}$ is 
convex in $\Delta(\alpha,\mu r)$.  Since this is true for each $\alpha\in\Bbb D$, the conclusion 
is that $h$ is uniformly locally convex.
\qed\enddemo
\demo{Proof of Theorem 3 }  We showed in [4] that  $\|{\Cal S}f\|<\infty$ if and only if 
$\|{\Cal S}h\|<\infty$.  Therefore, if $\|{\Cal S}f\|<\infty$, then $\|{\Cal S}h\|<\infty$, and so $h$ 
is uniformly locally univalent, by the theorem of B. Schwarz.  In other words, $h$ is univalent 
on every pseudohyperbolic disk $\Delta(\alpha,r)$ for some fixed radius $r$.  Then by the 
classical result on radius of convexity (see [8], p. 44), $h$ maps every disk 
$\Delta(\alpha,(2-\sqrt{3})r)$ to a convex domain.  It now follows from Theorem B that $f$ is 
univalent in each disk $\Delta(\alpha,(2-\sqrt{3})r)$.  Thus $f$ is uniformly locally univalent.  

     Conversely, suppose the  harmonic mapping $f$ is uniformly locally univalent in $\Bbb D$.  
Then by Theorem C its analytic part $h$ is uniformly locally convex, hence uniformly locally 
univalent.  Therefore, $\|{\Cal S}h\|<\infty$ by Kraus' theorem, as discussed at the beginning 
of this section.  It now follows from our result in [4] that $\|{\Cal S}f\|<\infty$.
\qed\enddemo

     As a corollary of the proof, we are able to establish a numerical bound on $\|{\Cal S}f\|$ 
for univalent harmonic mappings $f$, analogous to Kraus' bound $\|{\Cal S}f\|\leq6$ for 
analytic univalent functions in the disk.  By M\"obius invariance we may assume without loss of generality that the harmonic mapping $f=h+\overline{g}$ belongs to the class $S_H$.   Then 
as shown in the proof of Theorem C, its analytic part $h(z)$ is convex in the disk $|z|<\mu$, 
where $\mu=\lambda - \sqrt{\lambda^2 - 1}$.  Thus the function $H(z)=h(\mu z)$ is convex in 
$\Bbb D$, so it has Schwarzian norm $\|{\Cal S}H\|\leq2$, by a result of Nehari [15].  Since 
$\|{\Cal S}H\|=\mu^2\|{\Cal S}h\|$, it follows that $\|{\Cal S}h\|\leq2/\mu^2$.  Consequently, 
the estimate $\lambda<49$ shows that $\|{\Cal S}h\|< 19,204$.  On the other hand, 
a result of Pommerenke [18] implies that 
$$
\|{\Cal S}f\| \leq  \|{\Cal S}h\| + 2\left(1 + \tfrac12 \|{\Cal S}h\| \right)^{1/2} + 7\,, 
$$
as we showed in [4].  Inserting the preceding estimate $\|{\Cal S}h\| < 19,204$, we obtain the 
absolute bound $\|{\Cal S}f\| < 19,407$ for all harmonic mappings $f$ that are univalent in 
$\Bbb D$ and have dilatation that is a perfect square.
   
      It is an open problem to determine the sharp bound.  We showed in [4] that 
$\|{\Cal S}f\|\leq45$ for all mappings $f$ with dilatation $\omega=q^2$ that are convex 
in the horizontal direction.  We also observed that the horizontal shear of the Koebe function 
with dilatation $\omega(z)=z^2$ has Schwarzian 
$$ 
{\Cal S}f = - \,4\left(\frac{1}{1-z} + \frac{\overline{z}}{1 + |z|^2} \right)^2\,,
$$
from which an easy calculation gives $\|{\Cal S}f\|=16$.  These results are unchanged if the 
Koebe function is sheared with dilatation $\omega(z)=e^{i\theta}z^2$ for any $\theta$.  
Therefore, since the Koebe function maximizes the Schwarzian norm for analytic univalent 
functions, it is reasonable to conjecture that $\|{\Cal S}f\|\leq16$ for all univalent harmonic 
mappings in the disk whose dilatation is a perfect square.  

\bigpagebreak
\flushpar
{\bf \S 4.  Bounds on valence of harmonic lifts.}   
\smallpagebreak

     Theorems 1 and 2 extend readily to the lifts of harmonic mappings to minimal surfaces.  
In [3] we obtained the following generalization of Nehari's theorem.  

\proclaim{Theorem D}  Let $f=h+\overline{g}$ be a harmonic mapping of the 
unit disk, with conformal parameter $e^{\sigma(z)}=|h'(z)|+|g'(z)|\neq0$ and dilatation 
$g'/h'=q^2$ for some meromorphic function $q$.  Let $\widetilde{f}$ denote 
the Weierstrass--Enneper lift of $f$ to a minimal surface with Gauss curvature 
$K=K(\widetilde{f}(z))$ at the point $\widetilde{f}(z)$.  Suppose that 
$$
|{\Cal S}f(z)| + e^{2\sigma(z)} |K(\widetilde{f}(z))| \leq 2p(|z|)\,, \qquad z\in\Bbb D\,,
$$
for some Nehari function $p$. Then $\widetilde{f}$ is univalent in $\Bbb D$.
\endproclaim

     The valence estimates for analytic functions in Theorems 1 and 2 have corresponding  
generalizations to harmonic lifts.  

\proclaim{Theorem $1'$}  Let $f=h+\overline{g}$ be a harmonic mapping of the unit disk 
with conformal parameter $e^{\sigma(z)}=|h'(z)|+|g'(z)|\neq0$ and dilatation 
$g'/h'=q^2$ for some meromorphic function $q$, and let $\widetilde{f}$ be its lift  
to a minimal surface with Gauss curvature $K$.   Suppose that 
$$
|{\Cal S}f(z)| + e^{2\sigma(z)} |K(\widetilde{f}(z))| \leq C\,, \qquad z\in\Bbb D\,, 
$$
for some constant $C>{\pi}^2/2$. Then $|\alpha-\beta|\geq\sqrt{2/C}\,\pi$ for any pair 
of points $\alpha,\beta\in{\Bbb D}$ where $\widetilde{f}(\alpha)=\widetilde{f}(\beta)$.  
Consequently, the lift  $\widetilde{f}$ has finite valence and meets any given point 
at most $\left(1 + \frac{\sqrt{2C}}{\pi}\right)^2$ times.  
\endproclaim

\proclaim{Theorem $2'$}  Let a harmonic mapping $f=h+\overline{g}$ be as in Theorem $1'$ 
but satisfy the inequality  
$$
|{\Cal S}f(z)| + e^{2\sigma(z)} |K(\widetilde{f}(z))| \leq \frac{2C}{1-|z|^2}\,, \qquad z\in\Bbb D\,, 
$$
for some constant $C>2$. Then its lift  $\widetilde{f}$ has finite valence $N=N(C)\leq AC\log C$, 
where $A$ is some absolute constant. 
\endproclaim

     The proofs of Theorems $1'$ and $2'$ reduce ultimately to the same consideration of 
zeros of solutions to differential equations as in the proofs of Theorems 1 and 2.  Here the 
link with differential equations and the Sturm theory comes from a result of Chuaqui and 
Gevirtz [5], as developed in our earlier work [3,4].  The details are relatively straightforward 
and will not be pursued here.

\enddocument